\documentclass{amsart}
\usepackage{amssymb}
\usepackage{amscd}
\usepackage{bbm}

\newcommand\Rs{{\mathbb R}}
\newcommand\End{\mathop {\fam 0 End}\nolimits}

\newcommand\id{\mathbbm{1}}

\newcommand\Exp{\mathop {\fam 0 Exp}\nolimits}

\begin{document}

{\large

\title[DIFFERENTIAL OF THE EXPONENTIAL MAP]
{THE TAYLOR SERIES RELATED TO THE DIFFERENTIAL OF THE EXPONENTIAL MAP}
%\thanks{Supported in part by the grant RFFI-08-01-92001}
\author{Alexey~V.~Gavrilov}

\maketitle

\begin{abstract}
 In this paper we study the Taylor series of an
operator-valued function related to the differential of the
exponential map. For a smooth manifold $\mathcal{M}$ with a
torsion-free affine connection the operator $\mathcal{E}_p(v)$
acting on the space $T_p\mathcal{M}$ is defined to be the
composition of the differential of the exponential map at $v\in
T_p\mathcal{M}$ with parallel transport to $p$ along the geodesic.
The Taylor series of $\mathcal{E}_p$ as a function of $v$ is found
explicitly in terms of the curvature tensor and its high order
covariant derivatives at $p$. \end{abstract}

{{\it Key words:} affine connection, exponential map, Taylor
series.}

{{\it 2010 Mathematical Subject Classification:} 53B05.}

\section{INTRODUCTION}

Let $(\mathcal{M},\nabla)$ be a smooth manifold with a torsion-free
affine connection. For a point $p\in\mathcal{M}$, consider the
exponential map
$$\Exp_p:T_p\mathcal{M}\to\mathcal{M}.$$
Let $v\in T_p\mathcal{M}$ be a vector at a fixed point $p$ and
$\gamma_v$ be the corresponding geodesic $\gamma_v(t)=\Exp_p(tv)$.
 As $T_p\mathcal{M}$ is a linear space, the
differential of the exponential map may be considered a linear map
of the form
$$d_v\Exp_p:T_p\mathcal{M}\to T_{\Exp_pv}\mathcal{M}.$$

Denote by
$$\mathcal{I}_p(v):T_p\mathcal{M}\to T_{\Exp_pv}\mathcal{M}$$
the operator of parallel transport along the geodesic $\gamma_v$.
Let $\mathcal{E}_p$ be the map
$$\mathcal{E}_p:T_p\mathcal{M}\to \End_{\Rs}(T_p\mathcal{M}),$$
defined by
$$\mathcal{E}_p(v)=\mathcal{I}_p(v)^{-1}\circ d_v\Exp_p.$$
This is a smooth map  between two linear spaces, hence its Taylor
series is well defined. Our goal is to find this series explicitly
in terms of the curvature tensor.

For the best of author's knowledge, the problem has never been
considered in this generality. However, for a symmetric space
$\mathcal{M}$ the answer is well known \cite{H1} (see also
\cite{H3}, Ch. IV, Theorem 4.1). In our terms
$$\mathcal{E}_p(v)=\sum_{k=0}^{\infty}\frac{1}{(2k+1)!}r_0^k(v);\eqno{(1)}$$
the operator $r_0$ is defined by
$$r_0(v):w\mapsto R_p(v,w)v,\, w\in T_p\mathcal{M},$$
where $R$ is the curvature tensor. Note that this is a very special
case, because for a more general manifold the series depends not
only on the curvature tensor itself, but on its high order covariant
derivatives at the point $p$ as well.

Apparently, the only relevant result known for a general affine (or
Riemannian) manifold is the Helgason's formula \cite{H2} (in fact,
the formula is proved for an {\it analytic} manifold with an affine
connection). The formula is
$$d_v\Exp_pw=\frac{1-e^{-{\rm ad}v^*}}{{\rm ad}v^*}w^*\bigg{|}_{\Exp_pv},$$
where the adjoint refers to the Lie algebra of smooth vector fields
and $v^*$ denotes the vector field defined by the condition
$$v^*_{\Exp_pu}=\mathcal{I}_p(u)v,\,u\in T_p\mathcal{M}$$
(the same for $w^*$). This, of course, is not a Taylor series, and
its relation to the problem remains quite obscure.

The author might also mention his own paper \cite{G}, where an
algorithm of computing the Taylor series of the inverse operator
$\mathcal{E}_p(v)^{-1}$ (denoted there by $H(v)$) is proposed. This
algorithm, however, is quite involved in comparison with the
explicit formula presented here.

{\bf Remark}. In one respect the Helgason's formula is more general
than the results obtained below: it is not necessary to assume the
connection to be torsion-free. It should be noted, however, that the
differential of the exponential map (as the map itself)  does not
depend on the torsion part of the affine connection. So, it is
natural to restrict ourselves to the torsion-free case. In a more
general case the Taylor series contains terms which depend on the
torsion tensor, but these terms actually have nothing to do with the
exponential map.

\section{THE SERIES}

Before we can formulate the main result we have to introduce some
notation. For $n\ge 0$ and $v\in T_p\mathcal{M}$ denote by
$R_{p,v}^{(n)}$ the $n$-th order covariant derivation of the
curvature tensor at the point $p$ in the direction $v$. That is,
$$R_{p,v}^{(n)}=v^{\otimes n}\cdot\left(\nabla^n R\right)_p.$$
(The common notation $\nabla^n_{v,\dots,v}$ for hight order
covariant derivation in the direction $v$ is not very convenient
when $n$ is variable. For this reason we use the contraction
notation. The sign `$\cdot$' on the right hand side denotes the
contraction of the polyvector $v^{\otimes n}$ with the tensor
$\nabla^n R$).

For $p,v,n$ as above denote by
$$r_{p,n}(v)\in \End_{\Rs}(T_p\mathcal{M})$$
the linear operator defined by
$$r_{p,n}(v):w\mapsto R_{p,v}^{(n)}(v,w)v,\,w\in T_p\mathcal{M}.$$
For the sake of convenience we usually omit the point and,
sometimes, the vector:
$$r_n=r_n(v)=r_{p,n}(v).$$
For example,
$$r_0:w\mapsto R_p(v,w)v,\,r_1:w\mapsto(\nabla_v R)_p(v,w)v,$$
etc. Obviously, $r_n$ as a function of $v$ is homogeneous of degree
$n+2$
$$r_n(tv)=t^{n+2}r_n(v),\,t\in\Rs.$$

We also need an appropriate notation for compositions of operators
of this kind. Call a finite sequence of nonnegative integers a {\it
list}. The set of all lists (including the empty one) is denoted by
$\Lambda$. The empty list is denoted by the symbol $\varnothing$; to
write down a nonempty list we use square brackets; for example,
$$\nu=[2,0,1]\in\Lambda.$$
For every list $\nu\in\Lambda$ there is a corresponding operator
$r_{\nu}\in \End_{\Rs}(T_p\mathcal{M})$ which is a composition of
simple operators $r_n$. Namely, if $\nu=[n_1,n_2,\dots,n_k]$ then
$$r_{\nu}=r_{n_1}r_{n_2}\dots r_{n_k}.$$
By definition, $r_{\varnothing}=\id$.

We shall need three number functions on lists: the factorial $\nu!$,
the degree $|\nu|$ and the denominator $c_{\nu}$. By definition,
$\varnothing!=1,|\varnothing|=0$. For $\nu=[n_1,n_2,\dots,n_k]$ we
define the factorial and the degree as follows
$$\nu!=\prod_{j=1}^k (n_j!),\,|\nu|= 2k+\sum_{j=1}^k n_j.$$
Obviously,
$$r_{\nu}(tv)=t^{|\nu|}r_{\nu}(v),\,t\in\Rs,$$
which is where the term ``degree'' comes from.

The denominator is defined by $c_{\varnothing}=1$ and a recurrent
relation
$$c_{\nu}=|\nu|(|\nu|+1)c_{\nu^{\prime}},\,\nu\in\Lambda\setminus\{\varnothing\},$$
where $\nu^{\prime}$ is obtained from $\nu$ by omitting the first
element from the list. For example, $|[2,0,1]|=2\cdot 3+2+0+1=9$,
hence
$$c_{[201]}=9\cdot 10c_{[01]}=90\cdot 30c_{[1]}=
90\cdot 30\cdot 12c_{\varnothing}=32400.$$

Now we are able to formulate the main result.

{\bf Theorem} {\it Let $(\mathcal{M},\nabla)$ be a smooth manifold
with a torsion-free affine connection. Let $p\in\mathcal{M}$ and
$v\in T_p\mathcal{M}$. Then for any $n\ge 0$ the following equality
holds
$$\frac{1}{n!}\frac{d^n}{dt^n}\mathcal{E}_p(tv)\bigg{|}_{t=0}=
\sum_{|\nu|=n}\frac{1}{\nu!c_{\nu}}r_{\nu}(v), \eqno{(2)}$$ where
the sum on the right hand side is taken over all the lists $\nu\in
\Lambda$ of degree $n$.}

In other words, we have the Taylor series
$$\mathcal{E}_p(v)=\sum_{\nu\in\Lambda}\frac{1}{\nu!c_{\nu}}r_{\nu}(v). \eqno{(3)}.$$
For example, there are 13 lists of degree not greater than 6, namely
$\varnothing,[0],[1],[2],$
$[0,0],[3],[1,0],[0,1],[4],[2,0],[1,1],[0,2],[0,0,0]$. Computing the
corresponding coefficients, we have the

{\bf Corollary } {\it The function $\mathcal{E}_p$ can be expressed
as follows
$$\mathcal{E}_p(v)=\id+\frac{1}{6}r_0+\frac{1}{12}r_1+\frac{1}{40}r_2+\frac{1}{120}r_0^2+
\frac{1}{180}r_3+\frac{1}{180}r_1r_0+\frac{1}{360}r_0r_1+\frac{1}{1008}r_4+$$
$$+\frac{1}{504}r_2r_0+\frac{1}{504}r_1^2+\frac{1}{1680}r_0r_2+
\frac{1}{5040}r_0^3+\rho_7(v),$$ where $\rho_7(v)=O(|v|^7)$ as $v\to
0.$ }

Note that for a list $\nu=[0,0,\dots,0]$ which consists of $k$ zeros
we have $|\nu|=2k$, hence
$$c_{\nu}=2k(2k+1)c_{\nu^{\prime}}=(2k+1)!.$$
Thus, the coefficient at the term $r_{\nu}=r_0^k$ in the Taylor
series is equal to $1/(2k+1)!$, in agreement with (1).

It may be shown that there are no nontrivial algebraic relations
between the operators $r_n$ on a general manifold. So, the series
(3) is unique as a formal series in non-commutative variables $r_n$.

\section{PARALLEL TRANSPORT OF THE CURVATURE TENSOR}

We shall need some properties of the covariant derivation along a
geodesic. For our purpose it is convenient to consider a connection
on a vector bundle instead of an affine connection (which is a
connection on the tangent bundle). Let $E\to \mathcal{M}$ be a
smooth vector bundle on a smooth manifold. The space of smooth
sections of $E$ is denoted by $C^{\infty}(\mathcal{M},E)$. The
bundle is provided with a connection $\nabla$, which is a linear map
$$\nabla:C^{\infty}(\mathcal{M},E)\to C^{\infty}(\mathcal{M},T^*\mathcal{M}\otimes E),$$
satisfying the Leibniz rule
$$\nabla fu=df\otimes u+f\nabla u,\,f\in C^{\infty}(\mathcal{M}),
u\in C^{\infty}(\mathcal{M},E).$$ Note that the manifold itself is
not supposed to have a connection for now.

Let $\gamma:I\to \mathcal{M}$ be a smooth curve, where $I\subset\Rs$
is an interval. The induced connection on the restricted bundle
$\gamma^*E\to I$ is usually denoted by the sign $D$.  It is
convenient to use this connection in the form of a differential
operator $\frac{D}{dt}$, where $t\in I$ is the parameter on the
curve. This operator is  called a {\it covariant derivation along
the curve} $\gamma$. In other words, for any section $u\in
C^{\infty}(\mathcal{M},E)$ we have the equality
$$\frac{D}{dt}\gamma^*u=\dot{\gamma}\cdot\gamma^*\nabla u, \eqno{(4)}$$
where
$$\dot{\gamma}(t)=\frac{d}{dt}\gamma(t)\in T_{\gamma(t)}\mathcal{M},\,t\in I.$$
To be less pedantic, one may write (4) as an operator identity
$$\frac{D}{dt}=\nabla_{\dot{\gamma}}.$$

We shall need the following simple, but useful

{\bf Lemma 1} {\it Let $(\mathcal{M},\nabla)$ be a smooth manifold
with an affine connection and $E\to \mathcal{M}$ be a smooth vector
bundle with a connection. Let $\gamma:I\to \mathcal{M}$ be a
geodesic. Then for any $n\ge 1$ and any smooth section $u$  of the
bundle $E$
$$\frac{D^n}{dt^n}\gamma^*u=\dot{\gamma}^{\otimes n}\cdot\gamma^*\nabla^nu,\eqno{(5)}$$
where $\frac{D}{dt}$ is the covariant derivation along $\gamma$.}

The equality (5) may also be written in the operator form
$$\frac{D^n}{dt^n}=\dot{\gamma}^{\otimes n}\cdot\nabla^n.$$
Apparently, this formula is well known, but the author does not know
a proper reference. For this reason, we present here a proof. Note
that for $n>1$ the section $\nabla^nu$ on the right hand side of (5)
depends on the affine connection on the manifold while the left hand
side depends on the curve $\gamma$ only.

{\it Proof.} It is an almost trivial fact that the derivation
$\frac{D}{dt}$ can be canonically extended to the products of $E$
with tensor bundles on $\mathcal{M}$ and inherits all the common
properties of the covariant derivation. In particular, it is
compatible with contraction. For $n=1$ the formula (5) is the
definition of $\frac{D}{dt}$. If it is valid for some $n$, then we
have
$$\frac{D^{n+1}}{dt^{n+1}}\gamma^*u=\frac{D}{dt}\frac{D^n}{dt^n}\gamma^*u=
\frac{D}{dt}\dot{\gamma}^{\otimes n}\cdot\gamma^*\nabla^nu= \left(
\frac{D}{dt}\dot{\gamma}^{\otimes n}\right)\cdot\gamma^*\nabla^nu+
\dot{\gamma}^{\otimes n}\cdot\frac{D}{dt}\gamma^*\nabla^nu.$$ By
assumption, $\gamma$ is a geodesic, hence
$$\frac{D}{dt}\dot{\gamma}=0$$
and the first term on the right hand side vanishes. Applying (4) to
the section $\nabla^nu\in
C^{\infty}(\mathcal{M},T^*\mathcal{M}^{\otimes n}\otimes E),$ we
have
$$\frac{D}{dt}\gamma^*\nabla^nu=\dot{\gamma}\cdot\gamma^*\nabla^{n+1}u.$$
Thus,
$$\frac{D^{n+1}}{dt^{n+1}}\gamma^*u=
\dot{\gamma}^{\otimes n+1}\cdot\gamma^*\nabla^{n+1}u,$$ and the
equality (5) follows by induction. \qed

The covariant derivation is closely related to parallel transport.
As above, denote by $\mathcal{I}_p=\mathcal{I}_{p,E}$ the operator
of parallel transport of the form
$$\mathcal{I}_p(v):E_p\to E_{\Exp_pv}$$
along the geodesic $\gamma_v(t)=\Exp_p(tv)$. By the definition of
parallel transport,
$$\frac{D}{dt}\mathcal{I}_p(tv)z=0,$$
where $\frac{D}{dt}$ is the covariant derivation along $\gamma_v$
and $z\in E_p$ is a constant. In a more general case, when $z=z(t)$
depends on the parameter $t$, we have the operator equality
$$\frac{D}{dt}\circ \mathcal{I}_p(tv)=\mathcal{I}_p(tv)\circ\frac{d}{dt}, \eqno{(6)}$$
which will be of use below.

After the above preliminaries we have come to the matter. Consider
the operator ${\mathcal R}_p(v)\in\End_{\Rs}(T_p\mathcal{M})$,
defined by
$${\mathcal R}_p(v)w=\mathcal{I}_p(v)^{-1}R_{\Exp_pv}(\mathcal{I}_p(v)v,
\mathcal{I}_p(v)w)\mathcal{I}_p(v)v, w\in T_p\mathcal{M},
\eqno{(7)}$$ where $R_x$ denotes the curvature tensor at a point
$x\in\mathcal{M}$. As is well known, parallel transport is
compatible with tensor operations. Thus, this operator can also be
written in the form
$${\mathcal R}_p(v)w=(\mathcal{I}_p(v)^{-1} R_{\Exp_pv})(v,w)v.$$
In the latter equality the parallel transport operator
$\mathcal{I}_p(v)$ is applied to the curvature tensor, i.e. it acts
on the bundle $T^1_3\mathcal{M}$.

{\bf Lemma 2} {\it For $v\in T_p\mathcal{M}$ and $n\ge 0$,
$$\frac{d^{n}}{dt^{n}}{\mathcal R}_p(tv)\bigg{|}_{t=0}=
n(n-1)r_{n-2}(v).\eqno{(8)}$$ }

In fact, $r_{-1}$ and $r_{-2}$ are not defined, but it is convenient
to consider them arbitrary operators. So, for $n=0$ and $n=1$ the
right hand side of (8) is zero.

{\it Proof.} Consider the geodesic $\gamma=\gamma_v$. For $w\in
T_p\mathcal{M}$ we have
$${\mathcal R}_p(tv)w=t^2(\mathcal{I}_p(tv)^{-1} R_{\gamma(t)})(v,w)v,$$
hence
$$\frac{d^{n}}{dt^{n}}{\mathcal R}_p(tv)\bigg{|}_{t=0}w=
n(n-1)\left(\frac{d^{n-2}}{dt^{n-2}}\mathcal{I}_p(tv)^{-1}
R_{\gamma(t)}\right)\bigg{|}_{t=0}(v,w)v.$$ By (6),
$$\frac{d^{n-2}}{dt^{n-2}}\mathcal{I}_p(tv)^{-1}
R_{\gamma(t)}=\mathcal{I}_p(tv)^{-1}\frac{D^{n-2}}{dt^{n-2}}R_{\gamma(t)}.$$

By Lemma 1,
$$\frac{D^{n-2}}{dt^{n-2}}R_{\gamma(t)}\bigg{|}_{t=0}=v^{\otimes n-2}\cdot\left(
\nabla^{n-2}R\right)_p=R_{p,v}^{(n-2)}.$$ Taking into account the
equality $\mathcal{I}_p(0)=\id$ and the definition of $r_{n-2}$, we
have (8).\qed

\section{THE JACOBI FIELD}

For given vectors $v,w\in T_p\mathcal{M}$ consider a family of
geodesics
$$\gamma_{s}(t)=\Exp_p(tv+tsw),$$
parametrized by a real number $s$ taken from a neighbourhood of
zero. Denote by $\gamma=\gamma_0=\gamma_v$ the geodesic
corresponding to $s=0$. Let $J\in
 C^{\infty}(I,\gamma^*T\mathcal{M})$ be the vector field defined by
$$J_t=\frac{d}{ds}\gamma_{s}(t)\bigg{|}_{s=0}\in T_{\gamma(t)}\mathcal{M},$$
which is equivalent to
$$J_t=d_{tv}\Exp_p tw.$$
It is well known that a field of this kind is a Jacobi field
\cite{KN}. This means that it satisfies the differential equation
$$\frac{D^2}{dt^2}J_t=R_{\gamma(t)}(\dot{\gamma}(t), J_t)\dot{\gamma}(t).\eqno{(9)}$$

{\bf Lemma 3} {\it For any $v\in T_p\mathcal{M}$ the operator
$\mathcal{E}_p(tv)$ satisfies the differential equation
$$\left(t^2\frac{d^2}{dt^2}+2t\frac{d}{dt}\right)\mathcal{E}_p(tv)=
{\mathcal R}_p(tv)\mathcal{E}_p(tv).\eqno{(10)}$$ }

{\it Proof.} By the definition of $\mathcal{E}_p$,
$$\mathcal{E}_p(tv)tw=\mathcal{I}_p(tv)^{-1}\circ d_{tv}\Exp_p tw=\mathcal{I}_p(tv)^{-1}J_t.$$
Thus, by (6) and (9),
$$\left(t\frac{d^2}{dt^2}+2\frac{d}{dt}\right)\mathcal{E}_p(tv)w=
\frac{d^2}{dt^2}\mathcal{E}_p(tv)tw=
\mathcal{I}_p(tv)^{-1}\frac{D^2}{dt^2}J_t=\mathcal{I}_p(tv)^{-1}
R_{\gamma(t)}(\dot{\gamma}(t), J_t)\dot{\gamma}(t).$$

On the other hand, by the definition of the operator ${\mathcal R}$,
we have the equality
$${\mathcal R}_p(tv)\mathcal{E}_p(tv)w=t\mathcal{I}_p(tv)^{-1}
R_{\gamma(t)}(\dot{\gamma}(t), J_t)\dot{\gamma}(t).$$ Comparing
these equalities and taking into account that they are valid for any
choice of $w$, we have (10).\qed

{\it Proof of Theorem.} For $n\ge 0$ denote the left hand side of
(2) by
$$\mathcal{E}_n=\frac{1}{n!}\frac{d^n}{dt^n}\mathcal{E}_p(tv)\bigg{|}_{t=0}.$$
The equality
$$\mathcal{E}_n=\sum_{|\nu|=n}\frac{1}{\nu!c_{\nu}}r_{\nu}\eqno{(11)}$$
for $n=0$ and $n=1$ can be verified directly:
$$\mathcal{E}_0=r_{\varnothing}=\id,\,\mathcal{E}_1=0.$$

Taking the $n$-th order derivative of the both sides of  (10) at
$t=0$, we have by Lemma 2 the equality
$$n(n+1)n!\mathcal{E}_n=\sum_{k=0}^n\binom{n}{k}k(k-1)r_{k-2}\cdot(n-k)!\mathcal{E}_{n-k}.$$
For $n\ge 2$ the equality takes the form
$$\mathcal{E}_n=\frac{1}{n(n+1)}\sum_{m=0}^{n-2}\frac{1}{m!}r_m\mathcal{E}_{n-m-2}.\eqno{(12)}$$
By induction, we may assume that the formula (11) is valid for the
operators $\mathcal{E}_{n-m-2}$ on the right hand side of (12). We
have then the equality
$$\mathcal{E}_n=\sum_{m=0}^{n-2}\sum_{|\mu|=n-m-2}\frac{1}{m!\mu!n(n+1)c_{\mu}}r_mr_{\mu}.$$
Denote $\nu=[m,\mu]$ (then $\nu^{\prime}=\mu$). One can see that the
double sum in the latter equality can be replaced by a single sum
taken over the lists $\nu$ of degree $n$. Taking into account the
obvious equalities
$$\nu!=m!\mu!,\,c_{\nu}=n(n+1)c_{\mu},\,r_{\nu}=r_mr_{\mu},$$
we obtain (11). \qed

{Alexey~V.~Gavrilov, Department of Physics, Novosibirsk State
University, 2 Pirogov Street, Novosibirsk, 630090, Russia}
\begin{flushright}
\email{gavrilov19@gmail.com}
\end{flushright}

}


\begin{thebibliography}{}
%\bibitem{M}{M.P. do CARMO},{\it Riemannian Geometry}, Birkhauser,Boston ,1992.
\bibitem{G}{A.V.GAVRILOV}, {\it The Leibniz formula for the covariant derivative and some
applications}(in Russian), Mat. Tr. (Matematicheskie Trudy){\bf
13}(2010), 63-84.
\bibitem{H1}{S. HELGASON}, {\it On Riemannian curvature of homogeneous spaces}, Proc. Amer.
Math. Soc. {\bf 9}(1958), 831-838.
\bibitem{H2}{S. HELGASON}, {\it Some remarks on the exponential mapping for
an affine connection}, Math. Scand. {\bf 9}(1961), 129-146.
\bibitem{H3}{S. HELGASON}, {\it Differential Geometry, Lie Groups and Symmetric Spaces},
Graduate Studies in Mathematics. {\bf 34}, AMS, Providence, Rhode
Island, 1978.
\bibitem{KN}{S. KOBAYASHI, K. NOMIZU}, {\it Foundations of Differential Geometry}, II,
Interscience Publisher, New York , 1963.


\end{thebibliography}
\end{document}